\documentclass[11pt,a4paper]{article}

\usepackage{caption2}
\usepackage{CJK}
\usepackage{tabls}
\usepackage{bm}
\usepackage{multirow}
\usepackage{amssymb}
\usepackage{amsfonts}
\usepackage{mathrsfs}
\usepackage{empheq}
\usepackage{amsmath}
\usepackage{amsthm}
\usepackage{graphicx}
\usepackage{color}
\usepackage{indentfirst}
\usepackage{hyperref}
\usepackage{float}
\usepackage{cases}
\usepackage[titletoc]{appendix}

\usepackage[colorinlistoftodos,english]{todonotes}

\allowdisplaybreaks%

\def\f{\frac}
\def\pa{\partial}

\def\e{\eqref}

\def\lab{\label}

\def\i1n{i=1,\cdots,n}
\def\j1n{j=1,\cdots,n}
\def\ij1n{i,j=1,\cdots,n}

\def\R{\mathbb R}

\def\C{\mathbb C}

\def \i{\mathrm i}

 \numberwithin{equation}{section}
\theoremstyle{definition}

 \newtheorem{thm}{\indent Theorem}[section]
 \newtheorem{cor}{\indent Corollary}[section]
 \newtheorem{lem}{\indent Lemma}[section]

 \newtheorem{rem}{\indent Remark}[section]

\theoremstyle{definition}

\newcommand{\be}{\begin{equation}}
\newcommand{\ee}{\end{equation}}
\newcommand{\beq}{\begin{equation*}}
\newcommand{\eeq}{\end{equation*}}

\begin{document}

\begin{CJK*}{GB}{gbsn}
\title{\bf On the asymptotic stability of wave equations coupled by velocities of anti-symmetric type}

\author{ Yan Cui\thanks{School of Mathematics(Zhuhai), Sun Yat-sen University, P. R. China. E-mail: \texttt{cuiy27@mail.sysu.edu.cn}.} 
\quad  Zhiqiang Wang\thanks{School of Mathematical Sciences and Shanghai Key Laboratory for Contemporary Applied Mathematics,
Fudan University, Shanghai 200433, P. R. China.
E-mail: \texttt{wzq@fudan.edu.cn}. }
}

\date{}
\maketitle


\begin{abstract} 
In this paper, we study the asymptotic stability of two wave equations coupled by velocities of anti-symmetric type via only one damping. We adopt the frequency domain method to prove that the system with smooth initial data is logarithmically stable, provided that the coupling domain and the damping domain intersect each other.  Moreover, we show, by an example, that this geometric assumption of the intersection is necessary for 1-D case.
\end{abstract}


  \section{Introduction and Main Results} \lab{section1.1}

Let $\Omega\in \R^n$ be a bounded domain with smooth boundary $\pa \Omega$. 
We are interested in the asymptotic stability of the following system of two wave equations with Dirichlet boundary condition:
\be
    \label{wave}
        \left \{ 
            \begin{split}
               & y_{tt}-\sum_{j,k=1}^n(g^{jk}(x)y_{x_j})_{x_k}+\alpha(x)z_t+\beta(x)y_t=0 \quad &&\text{in} ~(0,+\infty)\times \Omega \\
               & z_{tt}-\sum_{j,k=1}^n(g^{jk}(x)z_{x_j})_{x_k}-\alpha(x)y_t=0         \quad &&\text{in}~ (0,+\infty)\times \Omega\\
               & y=z=0 \quad  &&  \text{on} \ (0,+\infty)\times \pa \Omega\\ 
               & (y(0),y_t(0),z(0),z_t(0))=(y^0,y^1,z^0,z^1) \quad  &&\text{in}~ \Omega.
            \end{split}
        \right.
    \ee
Here the coefficients of elliptic operator $g^{jk}(.)\in C^1(\overline{\Omega};\mathbb{R})$ satisfy
    \be
    \label{principle_coefficient}
        g^{jk}(x)=g^{kj}(x),\quad \forall x\in \overline{\Omega},~j,k=1,2,\cdots,n,
    \ee
and 
    \be
    \label{chapter four positive}
        \sum_{j,k=1}^n g^{jk}\xi^j\overline{\xi}^k\geq a  |\xi|^2,\quad \forall (x,\xi^1,\cdots,\xi^n)\in \overline{\Omega}\times  \mathbb{C}^n
    \ee
for some constant~$a >0$. 

We assume that the coupling coefficient $\alpha\in L^{\infty}(\Omega)$ and the  damping coefficient $\beta \in L^{\infty}(\Omega)$ are both nonnegative,
 and furthermore 
\be \label{nonzero}
 \omega_{\alpha} \triangleq  \{ x\in \Omega | \alpha(x) \neq 0\} \neq \emptyset \quad \text{and} \quad  \omega_{\beta}  \triangleq  \{ x\in \Omega | \beta(x) \neq 0\} \neq \emptyset.
\ee


It is classical to consider the system \eqref{wave} as the following Cauchy problem in space 
$\mathcal{H} \triangleq H_0^1(\Omega)\times L^2(\Omega) \times H_0^1(\Omega)\times L^2(\Omega)$: 
\be
    \label{abstra-sys}
        \left \{ 
            \begin{split}
               & \frac{dU}{dt} =\mathcal{A}U \\
               & U(0) =U_0 \triangleq (y^0,y^1,z^0,z^1) \in \mathcal{H}
            \end{split}
        \right.
    \ee
with $U= (y,u,z,v)$ and the linear operator $\mathcal{A}: \mathcal{D}(\mathcal{A})\subset \mathcal{H} \rightarrow \mathcal{H}$ is defined as 
    \be
    \label{DefineA}
         \left\{
            \begin{split}
            &\mathcal{A}U=\big(u,\sum_{j,k=1}^n(g^{jk}(x)y_{x_j})_{x_k}-\alpha(x)z_t-\beta(x)u_t, v, \sum_{j,k=1}^n(g^{jk}(x)z_{x_j})_{x_k} +\alpha(x)y_t\big),\\
                 &\mathcal{D}(\mathcal{A})
                =(H^2(\Omega)\cap H_0^1(\Omega))\times H^1_0(\Omega) \times (H^2(\Omega)\cap H_0^1(\Omega))\times H^1_0(\Omega),
            \end{split}
        \right.
    \ee
It is easy to know from the  theory of linear operator semigroup \cite{Pazy83}  that the system \eqref{abstra-sys} has a unique  solution $ U(t)= e^{t\mathcal{A}} U_0 $ in $C^0([0,+\infty), \mathcal{H})$.
Then, we can define the total energy of the system of \eqref{wave} by 
\be
    \label{Energy}
         \mathbb{E}(y,z)(t)=\f{1}{2}\int_\Omega\big(\sum_{j,k=1}^n g^{jk} y_{x_j} \overline{y}_{x_k}+ |y_t|^2 \big)dx+\f{1}{2}\int_\Omega\big(\sum_{j,k=1}^n g^{jk} z_{x_j} \overline{z}_{x_k}+ |z_t|^2 \big)dx,
    \ee
    which implies  immediately the  equivalence
    $$\mathbb{E}(y,z)(t)  \sim \|(y,y_t,z,z_t)(t,\cdot)\|^2_{\mathcal{H}}.$$
Obviously,  the total energy is non-increasing:
    \be\lab{energy-decay}
    \f{d}{dt} \mathbb{E}(y,z)(t)=-\int_\Omega \beta(x)|y_t|^2 dx\leq0, \quad  \forall t \geq 0.     \ee   
    
 The questions that we are interested in the following: 
    \begin{itemize}
    \item   Under what conditions on  $\alpha$ and $\beta$,  the system \eqref{wave} is asymptotically stable?
    \item   If the system \eqref{wave} is stable, what is the decay rate of the total energy $\mathbb{E}(y,z)(t)$  as $t \rightarrow +\infty$ ?
    \end{itemize}

 More precisely, the main result that we obtain is the following theorem.
\begin{thm}\lab{mainthm}
Assume that  \eqref{principle_coefficient}-\eqref{chapter four positive}-\eqref{nonzero} hold.  Assume furthermore  that 
there exist  a constant $\delta >0$ and a nonempty open subset $ \omega_\delta \subset \omega_{\alpha}\cap \omega_{\beta}  \subset \Omega $ such that
\be \label{geo1}
 \inf_{\omega_\delta}\alpha  \geq \delta \text{\quad and\quad } \inf_{\omega_\delta}\beta  \geq \delta.
\ee
 Then, there exists a constant $C>0$, such that for any initial data $(y^0,y^1,z^0,z^1)\in \mathcal{D(A)}$, the energy of solution to \eqref{wave} satisfies:
\be
    \lab{logstable}
    \mathbb{E}(y,z)(t)\leq  \f{C}{\ln(t+2)}\| (y_0,y_1,z_0,z_1)\|^2_ { \mathcal{D(A)}}, \quad \forall t \geq 0.
    \ee
Moreover,  the system \eqref{wave} is  strongly stable in $\mathcal{H}$, i.e.,  for any initial data $(y^0,y^1,z^0,z^1)\in \mathcal{H}$, 
\be \label{eng-decay}
\lim_{t\rightarrow +\infty}\mathbb{E}(y,z)(t) =0.
\ee  
\end{thm}


\begin{rem} In Theorem \ref{mainthm}, if we  assume instead that the damping and coupling coefficients $\alpha, \beta$ are both continuous on $\overline \Omega$, then the assumption \eqref{geo1} can be simply replaced by the following geometric condition:
  \be \label{geome}
 \omega_{\alpha}\cap \omega_{\beta} \neq \emptyset.
 \ee

 Independently, based on frequency domain method and multiplier method, Kassem-Mortada-Toufayli-Wehbe \cite{Chiraz2019Local} proved strong stability \eqref{eng-decay} when two waves propagate at different speed under the assumption \eqref{geome}.  
One can also refer to \cite[Theorem 2.1]{Fu12} for indirect stability results of other coupled wave system by displacements under the same geometric conditions. 
\end{rem}

\begin{rem}
We provide an example  in Section \ref{sect:ex}  to show  that  the geometric assumption $\omega_{\alpha}\cap \omega_{\beta}  \neq \emptyset$ is necessary in general, which is different from the situation with coupling by displacements. One can refer to the open problem raised in \cite[Remark 2.2]{Fu12}.  As a supplement,  we also refer to reader by \cite[section 5.2.1.3]{KMTW2020}, some numerical examples have been provided to show that for some initial data System \eqref{wave} seems also strongly stable when  $\omega_{\alpha}\cap \omega_{\beta}=\emptyset$. 
\end{rem}

\begin{rem} The result on logarithmical stability in Theorem \ref{mainthm} is sharp. Indeed, if $\alpha \equiv 0$,  the system \eqref{wave} is decoupled into a dissipative system for $(y,y_t)$ which is only logarithmically stable (see \cite{Lebeau96}) and a conservative one for $(z,z_t)$. Hence one can not expect a faster decay rate than the logarithmical one for the coupled system \eqref{wave} no matter what the coupling $\alpha$ is.
\end{rem}

\begin{rem} In the setting of Theorem \ref{mainthm},  similar stability results still hold for system \eqref{wave} with other types  of boundary conditions, for instance,  Robin conditions or mixed Dirichlet-Neumann conditions \cite{Fu12}.  However,  there are no such stability results for the system with Neumann conditions, since all the constant states are equilibrium of the system and will stay at the equilibrium all the time. 
\end{rem}

In order to prove  the logarithmic stability of the system \eqref{wave} with regular initial data in Theorem \ref{mainthm}, 
we adopt the frequency domain approach  to prove certain  spectral estimates of the infinitesimal generator $\mathcal{A}$ of the solution semigroup. One can refer to \cite{Burq98,Lebeau96} for the case of  single wave equation and \cite{Fu12} for the case of wave systems. 
\begin{thm}\lab{thm-res}
Suppose that  the assumptions of Theorem \ref{mainthm} hold. Then there exists a constant $C>0$ such that
\be\label{OO} 
 \mathcal{O}_C \triangleq \Big\{\gamma \in \C \,\Big |\, -\frac{e^{-C|\Im \gamma|}}{C} 
 \leq  \Re \gamma \leq 0 \Big\}  \cap   \Big\{\gamma \in \C \, |\, |\gamma | \geq \frac{2}{ C} \Big\}  \subset \rho(\mathcal{A})  
 \ee
and the following estimate holds
\be\lab{res-est-main}
\|(\mathcal{A}-\gamma I)^{-1}\|_{\mathcal{L}(\mathcal{H})}\leq Ce^{C|\Im \gamma|}, \quad \forall \gamma \in \mathcal{O}_C .
\ee
\end{thm}


Obviously, the energy decay given by  \eqref{eng-decay} implies directly the fact $ Sp(\mathcal{A})  \subset \{\gamma \in \mathbb{C} | \Re\gamma<0\}$ and in particular, the origin $O \in Sp(\mathcal{A})$.  Since $\rho(\mathcal{A}) $ is an open set, 
then the corollary follows from  \eqref{OO}, upon choosing $C$ large enough, as a byproduct of Theorem \ref{thm-res}. 

\begin{cor} Suppose that  the assumptions of Theorem \ref{mainthm} hold. Then there exists a constant $C>0$ such that
 \be\lab{spect}
 Sp(\mathcal{A}) \subset \{\gamma \in \mathbb{C} |\Re \gamma< -\frac{e^{-C|\Im \gamma|}}{C}\}.
\ee
\end{cor}

\begin{rem} The proof of Theorem \ref{thm-res} is based on Global Carleman estimates (see \cite{FLZ19}), which is quite elementary and allowed to address many stabilization problems for the system with lower order terms. Moreover, it can be used to obtain explicit bounds on some estimates of decay rate or constants costs in terms of the coefficients. 
Roughly speaking, \eqref{res-est-main} is equivalent to an observable estimate with constant cost like  $e^{C|\gamma|}$ for coupled elliptic system, which seems quite natural to adopt Global Carleman estimates to obtain these type estimates 
(see  Lemma \ref{lem-key} in section \ref{section:res-est} for more details). 
\end{rem}

\begin{rem} We should point out that we can not directly adopt the approach in this paper to obtain the logarithmic stability of System \eqref{wave} when two waves have different propagating speed. Roughly speaking, one key step in the proof of important Lemma \ref{lem-key} is using an easy fact 
that $\pa_s p \cdot [\pa_{ss}q+\pa_i (g^{ij}\pa_j q)]+\pa_s q \cdot [\pa_{ss}p+\pa_i (g^{ij}\pa_j p)]=\pa_s [\pa_s p\pa_s q+
p \pa_i ( g^{ij}\pa_{j}\pa_s q)]+\pa_i (g^{ij}\pa_s q \pa_j p)-\pa_i(g^{ij}p\pa_j \pa_s q)
$, which can be used to give an estimate that $L^2$ norm of the coupling term with force terms can control the $H^1$ energy. 
However, this fact is invalid for the case of two waves with different propagating speed. 
\end{rem}

\subsection{Previous results}

There are a lot of results about asymptotic  stability or stabilization of wave equations.  
Among them, 
Rauch-Taylor \cite{RauchTaylor74} and  Bardos-Lebeau-Rauch \cite{ BardosLebeauRauch92} pointed out  that, the single damped wave equation is exponentially stable
 if and (almost) only if the Geometric Control Condition (\emph{GCC}) is satisfied: 
\emph{There exists $ T>0$  such that every geodesic flow touches the support set of damping term 
 before $T$.}
If the damping acts on a small open set but the GCC is not satisfied, Lebeau\cite{Lebeau96} and Burq\cite{Burq98} proved that the wave equation is  logarithmically stable for regular initial data. There are also many results about polynomial stability of a single wave equation with special condition on the damping domain, cf. \cite{BurqHitrik07,Phung07,Leautaud14,LeautaudLerner16}. Recently Jin proved in \cite{J20} the damped wave equation on hyperbolic surface with constant curvature is exponentially stable even if the damping domain is arbitrarily small.

As for the case of coupled wave equations or other reversible equations,  indirect stability is an important issue both in mathematical theory or in engineering application.  
Indeed, it arises whenever it is impossible or too expensive to damp all the components of the state, and it is
hence important to study stabilization properties of coupled systems with a reduced number of feedbacks. For
finite dimensional systems, it is fully understood thanks to the Kalman rank condition. While in the case
of coupled partial differential equations, the situation is much more complicated. It depends not only on the algebraic structure of coupling but also the  geometric properties of the damping and coupling domain.

Alabau-Boussouira first studied indirect stability of a weekly coupled wave system where the coupling is through the displacements. In \cite{Alabau02}, she adopted multiplier method to obtain polynomial stability for wave system with anti-symmetric type coupling under stronger geometric conditions for both the coupling and damping terms. Moreover, she proved that this result was sharp for coupling with displacement. In \cite{AlabauLeautaud12}, the polynomial stability results for  coupled systems under an abstract framework (including wave-wave system, wave-plate system etc.) were obtained under the conditions that both coupling and damping are localized and satisfy the Piecewise Multipliers Geometric Conditions (PMGC,   see \cite{Liu97}). For 1-D case, a sharp decay rate of polynomial stability was obtained by Riesz basis method in \cite{LiuRao07}.  In \cite{Fu12}, Fu adopted global Carleman estimates and frequency domain method to prove that system with coupling by displacement of symmetric type was logarithmically stable with the assumption that coupling domain intersects the damping domain. 

 The above results concern only the weakly coupled system. In \cite{AWY17}, 
 Alabau-Wang-Yu studied the indirect stability for wave equations coupled by velocities with a general nonlinear damping. 
 By multiplier method, they obtained various types of stability results, including exponential stability, under strong geometric condition on the coupling  and damping domains. They also points out, for the first time, that it is more efficient way to transfer the energy in case of coupling  by velocities compared to the case with coupling by displacements. For 1-D case with constant coefficients, the sharp decay rate was explicitly given in \cite{CuiWang16}. In \cite{Klein17}, Klein computes the best exponent for the stabilization of wave equations on compact manifolds. The coefficient he obtains is therefore solution of some ODE system of matrices. Kassem-Mortada-Toufayli-Wehbe \cite{Chiraz2019Local} studied a system of two wave equations coupled by velocities with only one localized damping,  the waves propagate at different speed and the positivity and smallness assumptions of the continuous coupling coefficient can be removed. They obtained a strong stability result with the assumption that coupling domain intersects the damping domain. Moreover, assuming coupling and damping coefficient belong to $W^{1,\infty}(\Omega)$ and the intersection of coupling domain and damping domain holds PMGC, based on Frequency Domain method and multiplier method,  they established an exponential energy decay when the waves propagate at the same speed and polynomial energy decay when the waves propagate at the different speed. Recently,  the exponential energy decay result has been generalized by Gerbi-Kassem-Mortada-Wehbe in \cite{KMTW2020} to the case that the intersection of  the coupling and the damping domain holds GCC.

To the authors' knowledge, most known indirect stability results are obtained under the geometric conditions that the damping domain intersects the coupling domain.
Indeed, this guarantees effectively  the  energy transmission in higher space dimension. It is remarkable that Alabau-Boussouira and L\'{e}auteau \cite{AlabauLeautaud12} proved an indirect stability result in 1-D case where the damping domain and the coupling domain are two intervals which do not intersect.

\subsection{Main contribution and ideas}
  As already pointed out in \cite{AWY17}, the energy transition is more efficient through the first order coupling (by velocities) compared to zero order coupling (by displacements). This is natural since the first order coupling effect can be seen as a bounded perturbation to the system while the zero order coupling is a compact one. Nevertheless,  one can not expect a faster decay (than logarithmical one) of the whole system \eqref{wave} even if  a first order coupling appears, because  the energy of the single wave equation with damping localized in small domain only decays logarithmically. In this sense, the stability results are sharp. 

Not surprisingly, the indirect stability result is obtained by assuming essentially the damping and coupling domain intersect. However, we give an example to show that this geometric condition is necessary in general for the wave system coupled by velocities. This is quite different from the system coupled by displacements, see \cite{AlabauLeautaud12}.

 As for the  proof of the main theorems,  we adopt the frequency domain approach to 
 reduce the stability problem to an estimate on resolvent  which can be obtained by global  Carleman estimates of an elliptic equation as in \cite{Fu12}.  
 Different from the system in \cite{Fu12},  there are no zero order terms explicitly  in the system \eqref{wave}. 
 In order to derive the $L^2$ energy  of the solution, we then need to make fully use of the coupling structure together 
 with Poincar\'{e} Inequality under homogeneous Dirichlet boundary condition.

  \subsection{ Organization of the paper}
The paper is organized as follows. In Section
\ref{preliminaries}, we recall some basic facts about frequency domain method and global Carlemann estimates for an elliptic equation. 
 In Section \ref{section:proof}, we give the proofs of Theorem \ref{mainthm} and Theorem \ref{thm-res} as well as the technical Lemma \ref{lem-key},
 which is crucial to the proof of Theorem \ref{thm-res}.
Finally in Section \ref{sect:ex}, we give an example of  system \eqref{wave} with $\omega_{\alpha}\cap \omega_{\beta} = \emptyset$, which is indeed unstable.

\section{Preliminaries}\lab{preliminaries}
In this section, we briefly recall the frequency domain method and global Carleman estimates for an elliptic equation.
\subsection{Frequency domain method}
Thanks to classical semigroup theory, $\mathcal{A}$ generates a~$C_0$-semigroup operator $\{e^{t\mathcal{A}}\}_{t\geq 0}$ on $\mathcal{H}$ such that $S(t)U_0$. 
    It is well-known that the logarithmic stability of the system \eqref{wave} is equivalent to a resolvent estimate \cite{Burq98, Lebeau96}. More precisely, we have
    \begin{lem}\lab{lem:resolvent-est}
    Let $\mathcal{A}$ be defined by \eqref{DefineA}. If
   \be
    \label{res-est1}
        \|(\mathcal{A}-i\sigma I)^{-1}\|_{\mathcal{L(H)}}\leq C e^{C|\sigma|},\quad \forall   \sigma\in \R, |\sigma|>1
        \ee
then, there exists $C>0$ such that  for any $U_0\in \mathcal{D}(\mathcal{A}^2) \triangleq \{U\in \mathcal{H}| \mathcal{A}U\in \mathcal{D}(\mathcal{A})\},$
\be \label{2.2}
\|e^{t\mathcal{A}}U_0\|_{\mathcal{H}}\leq \Big(\f{C}{\ln(t+2)} \Big)^2\|U_0\|_{\mathcal{D}(\mathcal{A}^2)}, ~~~\forall t\geq 0
.\ee
    \end{lem}
Obviously, \eqref{res-est-main} in Theorem \ref{thm-res} implies the assumption \eqref{res-est1} in Lemma \ref{lem:resolvent-est}. Once Theorem \ref{thm-res} is proved,
 the logarithmical decay estimate \eqref{logstable} in Theorem  \ref{mainthm} can be easily obtained by Calderon-Lions interpolation theorem and \eqref{2.2}.

\subsection{Global Carlemann estimates}

To obtain resolvent estimates \eqref{res-est-main}, we need to introduce the global Carlemann estimates for elliptic equation, see  \cite{Fu09,Fu12,FI96}.

Let $\omega_0$ be an open set such that ~$\omega_0 \subset \subset \omega_{\delta} \subset \omega_{\alpha}\cap \omega_{\beta} $. There exists ~$\hat{\psi}\in C^2(\overline{\Omega};\R)$~ such that
    \be
    \label{eq3.1}
        \hat{\psi}>0~ \text{\ in} ~\Omega ,\quad  \hat{\psi}=0~\text{\ on}~\pa \Omega 
         \quad \text{and } \quad  |\nabla\hat{\psi}|>0~\text{\ in}~ \overline{\Omega\backslash \omega_0}.
    \ee
Next, we introduce some  weight functions:
    \be
    \label{Multiplier l}
       \theta=e^l, \quad l=\lambda\phi,\quad \phi=e^{\mu\psi}
    \ee
with 
 \be \label{psi}
   \psi=\psi(s,x)\triangleq \f{ \hat{\psi}(x)}{|| \hat{\psi}||_{L^{\infty}}}+b^2-s^2, \quad  s\in [-b,b], x\in \overline{\Omega}.
    \ee
 Here $b>0$ and $\lambda, \mu, s \in\R$ are all constants. 

Let's consider a single elliptic equation:
\be
\lab{elliptic-equ}
\begin{cases}\displaystyle
w_{ss}+\sum_{j,k=1}^n(g^{jk}w_{x_j})_{x_k}=f,  &\text{\quad in\ } (-b,b)\times \Omega,\\
w=0,  &\text{\quad on\ }  (-b,b) \times \pa \Omega,\\
w(\pm b,\cdot)=0, &\text{\quad in\ } \Omega,
\end{cases}
\ee
where $g^{jk}(.) \in C^1(\Omega; \R)$~satisfy \eqref{principle_coefficient}-\eqref{chapter four positive}. Then for every $f\in L^2((-b,b)\times\Omega)$, 
the elliptic system \eqref{elliptic-equ} has a unique solution $w\in H_0^1((-b,b)\times\Omega)$.
Then we have the following global Carleman estimates of solution.
    \begin{thm}[\cite{Fu12}]
    \label{global-carlemann-estimates}
      Let $b\in (1,2]$ and ~$\theta, \phi \in C^2([-b,b]\times \overline{\Omega}; \R) $ be ~defined by \eqref{Multiplier l}-\eqref{psi}.
    Then there exists ~$\mu_0>0$,~such that for any~$\mu\geq \mu _0$, there exist ~$C=C(\mu)>0$~ and ~$\lambda_0=\lambda_0(\mu)$~ such that for any~$f\in L^2((-b,b)\times\Omega)$, the solution $w$ to system \eqref{elliptic-equ} satisfies
    \be
    \label{eq3.5}
        \begin{split}
            \lambda &\mu^2\int_{-b}^b\int_\Omega \theta^2\phi(a |\nabla w|^2+ |w_s|^2+\lambda^2\mu^2\phi^2|w|^2)dxds \\
            &\leq C  \int_{-b}^b\int_\Omega \theta^2|f|^2dxds 
             + C\lambda\mu^2 \int_{-b}^b\int_{\omega_0} \theta^2\phi(|\nabla w |^2+|w_s|^2+\lambda^2\mu^2\phi^2|w|^2)dxds,
        \end{split}
    \ee
for all~$\lambda\geq \lambda_0(\mu)$.
    \end{thm}

    \section{Proofs of  main theorems }\lab{section:proof}
   In this section, we give the proofs of the main theorems, i.e., Theorem \ref{mainthm} and Theorem \ref{thm-res}.  
   First in Subsection \ref{section:res-est}, we prove Theorem \ref{thm-res}, particularly the  resolvent estimates \eqref{res-est},  
   based on some interpolation estimates on elliptic equations.
  Then  by  Theorem \ref{thm-res} and  {\color{red}Calderon}-Lions interpolation inequality,  we conclude Theorem \ref{mainthm} in Subsection \ref{subsec3.2}. 
  Finally in Subsection \ref{section:elliptic}, we prove Lemma \ref{lem-key} concerning an interpolation inequality of coupled elliptic equations, 
  which is crucial to the proof of Theorem \ref{thm-res}.

\subsection{Proof of Theorem \ref{thm-res}}\lab{section:res-est}
 Throughout this subsection, we use
$C=C(\Omega,\alpha,\beta)$ to denote generic positive constants which may vary
from line to line (unless otherwise stated).  
%

%

Let  $F=(f^0,f^1,g^0,g^1)\in \mathcal{H}$ and $U_0 =(y^0,y^1,z^0,z^1)\in \mathcal{D(A)}$ be such that 
    \be
    \label{eq5.1}
        (\mathcal{A}-\gamma I)U_0=F,
    \ee
where $\gamma\in \C$ and $\mathcal{A}$ be given by \eqref{DefineA}. %
Then \eqref{eq5.1} 
is equivalent to
    \beq
        \begin{cases}
            -\gamma y^0+y^1=f^0\quad\quad&\text{in}~ \Omega,\\ \displaystyle
            \sum_{j,k=1}^n(g^{jk}y^0_{x_j})_{x_k}-\alpha(x)z^1-(\beta(x)+\gamma)y^1=f^1\quad\quad&\text{in}~ \Omega,\\
            -\gamma y^0+y^1=g^0\quad\quad&\text{in}~ \Omega,\\ \displaystyle
            \sum_{j,k=1}^n(g^{jk}z^0_{x_j})_{x_k}-\gamma y^1+\alpha(x)y^1=g^1\quad\quad&\text{in}~ \Omega,\\
            y^0=z^0=0\quad \quad&\text{on} ~\pa \Omega.
        \end{cases}
    \eeq
or furthermore
    \be
    \label{eq5.3}
        \begin{cases}\displaystyle
           \sum_{j,k=1}^n(g^{jk}y^0_{x_j})_{x_k}-\gamma^2 y^0-\gamma \alpha(x)z^0-\gamma \beta(x)y^0= F^0 \triangleq (\beta(x)+\gamma)f^0+f^1 
               \quad&\text{in}~ \Omega, \\ \displaystyle
           \sum_{j,k=1}^n(g^{jk}z^0_{x_j})_{x_k}-\gamma^2 z^0+\gamma \alpha(x)z^0= F^1 \triangleq\gamma g^0+g^1 \quad&\text{in}~ \Omega, \\
           y^0=z^0=0\quad \quad&\text{on} ~\pa \Omega, \\
           y^1=f^0+\gamma y^0, z^1=g^0+\gamma z^0 \quad&\text{in}~ \Omega.
        \end{cases}
    \ee
In order to prove \eqref{res-est-main}, it suffices to  prove that there exists a constant $C>0$, such that
\be \label{3.4}
\|(y^0,y^1,z^0,z^1)\|_{\mathcal{H}}\leq C e^{C|\Im\gamma|} \|(f^0,f^1,g^0,g^1)\|_{\mathcal{H}}, \quad \forall  \gamma\in \mathcal{O}_C . 
\ee

For this purpose, we set
    \be
    \label{eq5.4}
        p(s,x)=e^{i\gamma s}y^0(x), \quad q(s,x)=e^{i\gamma s}z^0(x), \quad (s,x)\in (-2,2) \times \Omega.
    \ee
Then~$p$~and~$q$~satisfy the following coupled elliptic equation
    \be\displaystyle
    \label{eq5.5}
        \begin{cases}\displaystyle
            p_{ss}+\sum_{j,k=1}^n(g^{jk}p_{x_j})_{x_k}+i \alpha(x)q_s +i \beta(x)p_s= G^0 \triangleq F^0 e^{i\gamma s}, \quad&\text{in}~ (-2,2)\times\Omega, \\
            \displaystyle
            q_{ss}+\sum_{j,k=1}^n(g^{jk}q_{x_j})_{x_k}-i \alpha(x)p_s= G^1 \triangleq F^1 e^{i\gamma s}, \quad&\text{in}~ (-2,2)\times\Omega, \\
            p=q=0,\quad &\text{on} ~(-2,2)\times\pa\Omega.
        \end{cases}
    \ee

Note that there are no boundary conditions on $s= \pm 2$ in the above system \eqref{eq5.5}.  
We have the following lemma on interpolation estimate, while its proof  is left in Subsection \ref{section:elliptic}. 
    \begin{lem}
    \label{lem-key}
     Under the assumption of Theorem \ref{mainthm}, there exists a constant~$C>0$~ such that for any $\lambda>0$ big enough, the solution $(p,q)$ to \eqref{eq5.5} with form \eqref{eq5.4} satisfies
    \be
    \label{eq4.1}
        \begin{split}
            ||p||_{H^1(Y)}+||q||_{H^1(Y)}\leq &Ce^{C\lambda} \big (||G^0||_{L^2(X)}+||G^1||_{L^2(X)}+ ||p||_{H^1(-2,2;L^2(\omega_{\delta}))} \big)\\
            &+Ce^{-2\lambda}(||p||_{H^1(X)}+||q||_{H^1(X)} ).
        \end{split}
    \ee
    where
    \be
     X\triangleq(-2,2)\times\Omega,\quad Y\triangleq(-1,1)\times\Omega, \quad Z\triangleq(-2,2)\times \omega_\delta.
  \ee 
    \end{lem}

On the other hand,  by \eqref{eq5.4}, we have 
    \be\displaystyle
    \label{3.10}
        \begin{cases}
            ||y^0||_{H_0^1(\Omega)}+ ||z^0||_{H_0^1(\Omega)}\leq C e^{C|\Im\gamma|} (||p||_{H^1(-1,1;H_0^1(\Omega))}+||q||_{H^1(-1,1;H_0^1(\Omega))}),\\
            ||p||_{H^1(-2,2;H_0^1(\Omega))}+||q||_{H^1(-2,2;H_0^1(\Omega))}\leq C e^{C|\Im\gamma|}(|\gamma|+1) ( ||y^0||_{H_0^1(\Omega)}+ ||z^0||_{H_0^1(\Omega)}),\\
            ||p||_{H^1(-2,2;L^2(\omega_\delta))}\leq  C e^{C|\Im\gamma|}|\gamma| ||y^0||_{L^2(\omega_\delta)}
        \end{cases}
    \ee
for some constant $C>0$. Combining  \eqref{3.10}  and  \eqref{eq4.1},  we get 
    \be\displaystyle
    \label{eq5.6}
       ||y^0||_{H_0^1(\Omega)}+ ||z^0||_{H_0^1(\Omega)}\leq C e^{C|\Im\gamma|} (||f^0,f^1,g^0,g^1)\|_{\mathcal{H}}+||y^0||_{L^2(\omega_\delta )}).
    \ee

 Next, we turn to estimate $||y^0||_{L^2(\omega_\delta)}$.
 Let $\zeta\in C^2_0(\Omega;\R)$  be a cutoff function such that 
 \be
0\leq  \zeta(x) \leq 1\text{\ \, in\ } \Omega \text{\quad and\quad }\zeta(x) \equiv 1  \text{\ in\ }\omega_\delta \subset \omega_{\alpha}\cap \omega_{\beta}.
 \ee
 Multiplying ~$y$-equation in \eqref{eq5.3} by $2\zeta \overline{y}^0$ and integrating by parts on $\Omega$ yield that %
    \be\displaystyle
    \label{eq5.7}
        \begin{split}
            &\int_\Omega \Big (-\sum_{j,k=1}^n(g^{jk}y_{x_j}^0)_{x_k}+\gamma^2y^0+\gamma \alpha(x)z^0 +\gamma \beta(x)y^0 \Big )\cdot2\zeta\overline{y}^0 dx\\
            =&2\gamma^2\int_\Omega\zeta|y^0|^2dx +2\int_\Omega\zeta\sum_{j,k=1}^n g^{jk}y^0_{x_j}\overline{y}^0_{x_k}dx
            -\int_\Omega\sum_{j,k=1}^n(g^{jk}\zeta_{x_j})_{x_k}|y^0|^2dx
            \\&+ 2\gamma \int_\Omega \alpha(x) \zeta \overline{y}^0 z^0dx +2  \gamma \int_\Omega \beta(x)\zeta|y^0|^2dx
            \\= & \int_\Omega F^0\cdot2\zeta\overline{y}^0 dx.
        \end{split}
    \ee
 Similarly, multiplying $z$-equation in \eqref{eq5.3} by $2\zeta\overline{z}^0$ and integrating by parts on $\Omega$ yield that 
    \be\displaystyle 
    \label{eq5.8}
        \begin{split}
        &\int_\Omega(-\sum_{j,k=1}^n(g^{jk}z_{x_j}^0)_{x_k}+\gamma^2z^0-\gamma \alpha(x)y^0)\cdot2\zeta\overline{z}^0 dx\\
         =&2\gamma^2\int_\Omega\zeta|z^0|^2dx+2\int_\Omega\zeta\sum_{j,k=1}^n g^{jk}z^0_{x_j}\overline{z}^0_{x_k}dx
         -\int_\Omega\sum_{j,k=1}^n(g^{jk}\zeta_{x_j})_{x_k}|z^0|^2dx
         \\&-2\gamma\int_\Omega \alpha(x) \zeta y^0 \overline{z}^0dx
         \\=& \int_\Omega F^1\cdot2\zeta\overline{z}^0 dx.
        \end{split}
    \ee
 Note that 
    \beq
       \Im (\gamma \overline{y}^0 z^0 -\gamma y^0\overline{z}^0)= \Im [ \gamma  (\overline{y}^0 z^0 - y^0\overline{z}^0)] = \Re \gamma \cdot 2\Im (\overline{y}^0 z^0).
    \eeq
Adding \eqref{eq5.7} to \eqref{eq5.8}  and  taking the imaginary part result in
    \beq
    \begin{split}
      4\Re \gamma \Im \gamma \int_\Omega \zeta(|z^0|^2+|y^0|^2)dx
      +2\Im \gamma\int_{\Omega}\zeta\beta(x)|y^0|^2dx
      &   + 2 \Re \gamma \int_\Omega \alpha\zeta \Im (\overline{y}^0 z^0) dx\\
         &=   \int_\Omega   2\zeta \cdot \Im ( F^0 \overline{y}^0 +F^1 \overline{z}^0 )dx. 
   \end{split}
    \eeq
Then  it follows by Cauchy-Schwartz inequality and the definition of $\beta, F^0,F^1$ that
    \beq
    \begin{split}
   | \Im \gamma|\int_{\omega_\delta}\beta(x)|y^0|^2dx  
 \leq  &  | \Im \gamma|\int_{\Omega}\zeta\beta(x)|y^0|^2dx\\
 \leq & C|\Re\gamma|(|\Im\gamma| +1) (||y^0||^2_{L^2(\Omega)}+ ||z^0||^2_{L^2(\Omega)}) \\
& +C (\|F^0\|_{L^2(\Omega)}^2 +\|F^1\|^2_{L^2(\Omega)})\\
 \leq & C|\Re\gamma|(|\Im\gamma| +1) (||y^0||^2_{L^2(\Omega)}+ ||z^0||^2_{L^2(\Omega)}) \\
& +C(|\gamma|+1) \|(f^0,f^1,g^0,g^1)\|_{\mathcal{H}}^2.
   \end{split}
    \eeq
Thus by  \eqref{eq5.6}, we get
\be \label{key}
\begin{split}
| \Im \gamma|\int_{\omega_\delta}\beta(x)|y^0|^2dx\leq 
& C|\Re\gamma|(|\Im\gamma| +1)e^{C|\Im\gamma|} ||y^0||^2_{L^2(\omega_\delta)} \\
& +C(|\Re\gamma||\Im\gamma| +|\gamma|+1)e^{C|\Im\gamma|}\|(f^0,f^1,g^0,g^1)\|_{\mathcal{H}}^2.
\end{split}
\ee
By definition of $\mathcal{O}_C$, we take  $C>0$ large enough such that  
 \be \label{3.18}
 C |\Re\gamma|(|\Im\gamma| +1)e^{C|\Im\gamma|}\leq \frac{\delta  |\Im \gamma| }{2} \text{\quad and \quad}  \Im \gamma > \frac{1}{C}
 \ee
for all $\gamma\in \mathcal{O}_C $. Note also the fact $\beta(x) \geq \delta $ a.e. in $\omega_{\delta}$. Then it follows from \eqref{key} and \eqref{3.18} that
    \be
    \label{eq5.11}
        ||y^0||_{L^2(\omega_{\delta})} \leq C e^{C|\Im \gamma|}  \|(f^0,f^1,g^0,g^1)\|_{\mathcal{H}}
    \ee
for some $C>0$ large enough.    
Combining  \eqref{eq5.6} and \eqref{eq5.11} gives 
    \beq
        ||y^0||_{H_0^1(\Omega)}+||z^0||_{H_0^1(\Omega)}\leq C e^{C|\Im \gamma|}  \|(f^0,f^1,g^0,g^1)\|_{\mathcal{H}}.
    \eeq
Since ~$y^1=f^0+\gamma y^0,z^1=g^0+\gamma z^0$, we have also
    \beq
    \begin{split}
        ||y^1||_{L^2(\Omega)}+||z^1||_{L^2(\Omega)}
        &\leq ||f^0||_{H_0^1(\Omega)}+||g^0||_{H_0^1(\Omega)}+| \gamma|(||y^0||_{L^2(\Omega)}+||z^0||_{L^2(\Omega)})\\
        &\leq Ce^{C|\Im \gamma|} \|(f^0,f^1,g^0,g^1)\|_{\mathcal{H}}.
    \end{split}
    \eeq
Hence the desired estimate \eqref{3.4} indeed holds for all $\gamma\in \mathcal{O}_C$. 

 Consequently, $\mathcal{A}-\gamma I$ is a bijection from $\mathcal{D}(\mathcal{A})$ to $\mathcal{H}$
which satisfies  the resolvent estimate \eqref{res-est-main}.
We conclude the proof of Theorem \ref{thm-res}. 

\subsection{Proof of Theorem \ref{mainthm}}\label{subsec3.2}

As a corollary of  Theorem \ref{thm-res},  there exists $C>0$ such that
  \be
    \label{res-est}
        \|(\mathcal{A}-i\sigma I)^{-1}\|_{L^2(\mathcal{H})}\leq C e^{C|\sigma|},\quad \forall   \sigma\in \R, |\sigma|>1.
                \ee
Then by Lemma \ref{lem:resolvent-est}, we have  for  $U_0\in \mathcal{D}(\mathcal{A}^2)$ that 
\be\lab{xx66}
\|e^{t\mathcal{A}}U_0\|_{\mathcal{H}}
\leq (\f{C}{\ln(t+2)})^2 \|U_0\|_{\mathcal{D}(\mathcal{A}^2)}, \quad   \forall t\geq0 .
\ee
On the other hand,   the  contraction of the semigroup $e^{t\mathcal{A}}$ implies that 
\be\lab{3.40}
\|e^{t \mathcal{A}}U_0\|_{\mathcal{H}}\leq  \|U_0\|_{\mathcal{H}}, \quad \forall t\geq0.\ee
Note that $\mathcal{D(A)}$ is an interpolate space between $\mathcal{D}(\mathcal{A}^2)$ and 
$\mathcal{H}$. 
Combing \eqref{xx66}-\eqref{3.40} and using Calderon-Lions interpolation theorem  (see \cite[p. 38, Example 1 and p. 44, Proposition 8]{RS72}),
 we conclude for all $U_0 \in \mathcal{D}(\mathcal{A})$ that 
 \be\lab{3.41}
\|e^{t \mathcal{A}}U_0\|_{\mathcal{H}}
\leq  \f{C}{\ln(t+2)} \|U_0\|_{\mathcal{D(A)}},
\quad \forall  t\geq 0,
\ee
which is equivalent to the logarithmical decay estimate \eqref{logstable}.

Finally we  conclude by \eqref{3.41}  and density argument that the system   \eqref{wave} is strongly stable, i.e., for all $U_0 \in \mathcal{H}$, 
 \beq 
\lim_{t\rightarrow +\infty} \|e^{t \mathcal{A}}U_0\|_{\mathcal{H}}=0.
\eeq  
The proof of Theorem \ref{thm-res} is complete.

\subsection{Proof of Lemma \ref{lem-key}}\lab{section:elliptic}
In this subsection, we give the proof of Lemma \ref{lem-key} 
  which plays a key role in proving Theorem \ref{thm-res}.
The proof  is divided into 6 steps. 
In this subsection, we denote $C>0$ various constants independent of $\lambda$ which can be different from one line to another.

{\bf $\circ $ Step 1}.  We derive a weighted estimate \e{3.35} for $(p,q)$, the solution  of \eqref{eq5.5}.

Note that there are no boundary conditions on $p(\pm 2,\cdot), q(\pm 2,\cdot)$ in \eqref{eq5.5}.  
Let us introduce a cutoff function ~$\varphi=\varphi(s)\in C_0^3((-b,b); \R)$   (see for instance \cite{KimZhang08}) such that
\be    \label{eq4.3}
0\leq  \varphi(s) \leq 1 \text{\ \,in\ }   [-b,b]   \text{\quad and\quad }  \varphi(s) \equiv 1  \text{\ \,in\ } [-b_0,b_0].
 \ee     
where the constants $b_0, b$  are given by:
    \be
    \label{eq4.4}
        b\triangleq\sqrt{1+\f{1}{\mu} \ln(2+e^\mu)},\quad b_0\triangleq\sqrt{b^2-\f{1}{\mu} \ln \Big(\f{1+e^\mu}{e^\mu} \Big)},\quad \forall \mu >0.
    \ee
Obviously, if $\mu > \ln 2$, then 
 \be  \label{b-b0}
  1 <b_0 <b <2. \ee     
Let
    \be
    \label{eq4.5}
        \hat{p}(s,x)=\varphi(s) p(s,x),\quad \hat{q}(s,x)=\varphi(s) q(s,x), \quad (s,x) \in (-b,b) \times \Omega.
    \ee
then, we consider the elliptic equations that $\hat p, \hat q$ satisfy in $(-b,b) \times \Omega$: 
    \be
    \label{define4.3}
        \begin{cases} \displaystyle
            \hat{p}_{ss}+\sum_{j,k=1}^n(g^{jk}\hat{p}_{x_j})_{x_k}
                 = \hat G^0, 
                 & \quad \text{in }  (-b,b) \times \Omega, \\
            \displaystyle
            \hat{q}_{ss}+\sum_{j,k=1}^n(g^{jk}\hat{q}_{x_j})_{x_k}
                 =\hat G^1, 
                 &\quad \text{in }  (-b,b) \times \Omega, \\
                 \hat p= \hat q=0,  &\text{\quad on\ }  (-b,b) \times \pa \Omega\\
                \hat p(\pm b,\cdot)= \hat q(\pm b,\cdot)=0, &\text{\quad in\ } \Omega,
        \end{cases}
    \ee
   where 
   \be \label{G0G1}
   \begin{split}
  &\hat G^0 \triangleq \varphi_{ss}p+2\varphi_sp_s
  +\varphi G^0+i\alpha(x)(\varphi_s q-\hat{q}_s)+i\beta(x)(\varphi_s p-\hat{p}_s),\\
  &\hat G^1 \triangleq \varphi_{ss}q+2\varphi_sq_s+\varphi G^1-i\alpha(x)(\varphi_s p-\hat{p}_s).
 \end{split}
   \ee 
  By applying Theorem \ref{global-carlemann-estimates} to both $\hat p$ and $\hat q$,  there exists ~$\mu_0> \ln 2$ such that for any $\mu\geq \mu_0$,
   there exist ~$C=C(\mu)>0$~ and ~$\lambda_0=\lambda_0(\mu)$~ such that for any~$\lambda\geq \lambda_0(\mu)$, we have
    \be
    \label{3.33}
        \begin{split}
            \lambda &\mu^2\int_{-b}^b\int_\Omega \theta^2\phi\big(a |\nabla \hat{p}|^2+ |\hat{p}_s|^2+\lambda^2\mu^2\phi^2|\hat{p}|^2\big)dxds \\
            \leq C&  \int_{-b}^b\int_\Omega \theta^2  |\hat G^0|^2dxds 
            + C\lambda\mu^2 \int_{-b}^b\int_{\omega_0} \theta^2\phi(|\nabla \hat{p}|^2+ |\hat{p}_s|^2+\lambda^2\mu^2\phi^2|\hat{p}|^2\big)dxds
        \end{split}
    \ee
    \be
    \label{3.34}
        \begin{split}
        \lambda &\mu^2\int_{-b}^b\int_\Omega \theta^2\phi  \big(a |\nabla \hat{q}|^2+ |\hat{q}_s|^2
        		+\lambda^2\mu^2\phi^2|\hat{q}|^2  \big)dxds \\
            \leq C&  \int_{-b}^b\int_\Omega \theta^2  |\hat G^1|^2 dxds  
            + C\lambda\mu^2 \int_{-b}^b\int_{\omega_0} \theta^2\phi   ( | \nabla \hat{q}|^2+ |\hat{q}_s|^2+\lambda^2\mu^2\phi^2|\hat{q} |^2 ) dxds,
        \end{split}
    \ee
where $\hat G^0,\hat G^1$ are given by \eqref{G0G1} and $\omega_0 \subset \subset \omega_{\delta} \subset \omega_{\alpha} \cap \omega_{\beta}$.
Adding \eqref{3.33} to \e{3.34} gives
 \be   \label{3.35} 
 I_0 \leq I_1+I_2
 \ee
where 
   \be \small \label{I012}
        \begin{split}
            I_0 & \triangleq   \lambda \mu^2\int_{-b}^b\int_\Omega \theta^2\phi\big(a |\nabla \hat{p}|^2+ |\hat{p}_s|^2+\lambda^2\mu^2\phi^2|\hat{p}|^2
               +a |\nabla \hat{q}|^2+ |\hat{q}_s|^2+\lambda^2\mu^2\phi^2|\hat{q}|^2 \big)dxds \\
               I_1 & \triangleq  C  \int_{-b}^b\int_\Omega \theta^2  (|\hat G^{0}|^2+ |\hat G^1|^2)  dxds    \\
              I_2 & \triangleq   C\lambda\mu^2 \int_{-b}^b\int_{\omega_0} \theta^2\phi(|\nabla \hat{p}|^2+ |\hat{p}_s|^2+     \lambda^2\mu^2\phi^2|\hat{p}|^2 +| \nabla \hat{q}|^2+ |\hat{q}_s|^2+\lambda^2\mu^2\phi^2|\hat{q} |^2  \big)dxds.
        \end{split}
    \ee

{\bf $\circ $  Step 2.}  We estimate  $I_0$  in \eqref{I012} from below. 

By the choice of $\theta,l,\phi$ in \eqref{Multiplier l} and $b, b_0$ in \eqref{eq4.4}, we have 
\be\lab{the2}
\theta \geq e^{\lambda(2+e^{\mu})},\quad \forall |s|\leq 1, x\in \Omega.
\ee
Then 
 \be    \small  \label{I0-bound}
      \begin{split}
            I_0
          \geq & \lambda \mu^2\int_{-b_0}^{b_0}\int_\Omega \theta^2\phi \big(a |\nabla p|^2+ |p_s|^2+\lambda^2\mu^2\phi^2|p|^2 
                    + a |\nabla q|^2+ |q_s|^2+\lambda^2\mu^2\phi^2|q|^2\big)dxds  \\
               \geq & \lambda  e^{2\lambda (2+e^{\mu})}  C(\mu) \int_{-1}^{1}\int_\Omega 
               	\big( |\nabla p|^2+ |p_s|^2+|p|^2 + |\nabla q|^2+ |q_s|^2+|q|^2\big)dxds\\
               = & \lambda  e^{2\lambda (2+e^{\mu})}  C(\mu) (\|p\|^2_{H^1(Y)} +\|q\|^2_{H^1(Y)}).
  \end{split}
    \ee

  {\bf $\circ $ Step 3.}      
  We estimate  $I_1$ in \eqref{I012} from above.

Using Cauchy-Schwarz inequality, we get easily 
  \begin{align} \label{I1}
  I_1   \leq I_{11}+I_{12}+I_{13}  
  \end{align}
where
   \be \small
    \label{3.36}
        \begin{split}
	I_{11} & \triangleq  C  \int_{-b}^b\int_\Omega \theta^2 (|\hat{G}_1^0|^2+|\hat{G}_1^1|^2)dxds,\\
           &\quad \text{with \ }  \hat{G}_1^0=\varphi_{ss}p+2\varphi_sp_s+i\alpha(x)\varphi_s q  +i\beta(x)\varphi_s p,
           \  \hat{G}_1^1= \varphi_{ss}q+2\varphi_sq_s-i\alpha(x)\varphi_s p|^2, \\
          I_{12} & \triangleq  C  \int_{-b}^b\int_\Omega \theta^2  ( | i\alpha(x)\hat{q}_s + i\beta(x)\hat{p}_s |^2 +| i\alpha(x) \hat{p}_s|^2)dxds, \\
           I_{13} & \triangleq  C \int_{-b}^b\int_\Omega \theta^2 \varphi^2 (|G^0|^2  + |G^1|^2)  dxds 
                   \end{split}
    \ee 
Here we denote $\hat G_1^0, \hat G_1^1$ some terms concerning the derivatives of $\varphi$ in $\hat G^0, \hat G^1$, which are useful for the estimation below.    

By the choice of $\theta,l,\phi$ in \eqref{Multiplier l} and $b, b_0$ in \eqref{eq4.4}, we have 
\be\lab{the1}
\theta \leq e^{\lambda(1+e^{\mu})}, \quad \forall b_0\leq|s|\leq b, x\in \Omega.
\ee
Then
 \be   \displaystyle\small \label{I11}
      \begin{split}
            I_{11}=&C\int_{(-b,-b_0)\bigcup(b_0,b)}\int_{\Omega}\theta^2
            		 \big( |\varphi_{ss}p+2\varphi_sp_s+i\alpha(x)\varphi_s q  +i\beta(x)\varphi_s p|^2 \\
            &\hspace{36mm}+|\varphi_{ss}q+2\varphi_sq_s-i\alpha(x)\varphi_s p|^2 \big)  dxds \\
             \leq&  e^{2\lambda(1+e^{\mu})}C(\mu)\int_{(-b,-b_0)\bigcup(b_0,b)}\int_{\Omega} \big( |p |^2 +|p_s|^2+| q|^2+|q_s|^2\big)dxds           \\
              \leq& e^{2\lambda(1+e^{\mu})}C(\mu) (\|p\|^2_{H^1(X)} +\|q\|^2_{H^1(X)})  .
  \end{split}
    \ee

Obviously, 
\be  \label{I12}
I_{12} \leq   C  \int_{-b}^b\int_\Omega \theta^2  ( |\hat{p}_s|^2+ |\hat{q}_s|^2 )dxds  \leq \frac{C(\mu)}{\lambda} I_0\leq \lambda^{-\f 12} I_0,
\ee
which can be absorbed by $I_0$ if $\lambda$ is large enough.

Therefore, 
 \be   \displaystyle \label{I1-bound}
      \begin{split}
            I_{1}\leq & e^{2\lambda(1+e^{\mu})}C(\mu) (\|p\|^2_{H^1(X)} +\|q\|^2_{H^1(X)})
            +\lambda^{-\f 12}  I_0 \\
            &+ C \int_{-b}^b\int_\Omega \theta^2 \varphi^2 (|G^0|^2  + |G^1|^2)  dxds. 
  \end{split}
    \ee

{\bf $\circ $  Step 4.}  We estimate the localized term $I_2$ in \eqref{I012} from above.

We write it as 
 \be  
 \begin{split} \label{I2}
 I_2 = I_{21}+I_{22}+I_{23}+I_{24},
 \end{split}
 \ee
 where 
 \begin{align*}  
 I_{21} &\triangleq  C\lambda\mu^2 \int_{-b}^b\int_{\omega_0} \theta^2\phi(|\nabla \hat{p}|^2+ |\hat{p}_s|^2)  dxds,\quad
 I_{22} \triangleq  C \lambda^3\mu^4 \int_{-b}^b\int_{\omega_0} \theta^2 \phi^3 |\hat{p}|^2  dxds,\\
  I_{23} &\triangleq  C\lambda\mu^2 \int_{-b}^b\int_{\omega_0} \theta^2\phi(|\nabla \hat{q}|^2+ |\hat{q}_s|^2)  dxds,\quad
 I_{24} \triangleq  C \lambda^3\mu^4 \int_{-b}^b\int_{\omega_0} \theta^2 \phi^3 |\hat{q}|^2  dxds.
 \end{align*}
  
 Recall the definitions $\hat p, \hat q$ in \eqref{eq4.5}, $\varphi$ in  \eqref{eq4.3} and  $\theta,l,\phi$ in
 \eqref{eq3.1}-\eqref{psi}.  It's easy to check that
    \be
    \label{eq4.7}
        \begin{split}
           & l_s=-2\lambda\mu s\phi,\quad l_{x_j}=\lambda\mu\phi\psi_{x_j},\quad l_{x_js}=-2\lambda\mu^2 s\phi\psi_{x_j}, \\
           & l_{ss}=4\lambda\mu^2 s^2\phi{\color{blue} -}2\lambda\mu\phi,\quad l_{x_jx_k}=\lambda\mu^2\phi\psi_{x_j}\psi_{x_k}+\lambda\mu\phi\psi_{x_jx_k}.
        \end{split}
    \ee
Denote $\omega_k \ (k=1,2,3)$ some open subsets in ~$\Omega$ such that 
$\omega_0 \subset \subset \omega_1 \subset \subset \omega_2 \subset \subset\omega_3\subset \subset\omega_\delta$. 
Let ~$\eta_j\in C^3_0(\omega_j; \R)\ (j=1,2,3)$~ be suitable  cut-off functions such that
    \be
    \label{eq4.8}
            \eta_j(x) \equiv1\text{\  in\ } \omega_{j-1},\quad 
            0\leq \eta_j(x) \leq 1 \text{\  in\ }  \omega_j, \quad 
            \eta_j(x) \equiv 0  \text{\  in\ }  \Omega\backslash \omega_j.
                \ee
  Moreover, we choose further $\eta_2$ such that 
 \be \label{eta2}
 |(\eta_2)_{x_j}(\eta_2)_{x_k}|\leq C\eta_2, \quad \forall x \in \omega_2.
 \ee
for some constant $C>0$.
The existence of $\eta_2$ is shown at the end of the proof.

 {\bf $\circ $  Step 4.1.} 
 We  estimate 
   a weighted energy for $(\nabla\hat q, \hat q_s)$, i.e., 
   $\int_{-b}^b\int_{\omega_0}\theta^2\phi(|\nabla \hat{q}|^2+ |\hat{q}_s|^2)dxds$.
   
   By definition of $\eta_1$, it suffices to estimate 
   $ \int_{-b}^b\int_{\Omega} \theta^2\phi \eta_1^2(a |\nabla \hat{q}|^2+|\hat{q}_s|^2)dxds$.
%
To do this, we multiply $\hat q$-equation  in \eqref{define4.3} by  $\theta^2\phi\eta_1^2 \overline{\hat{q}}$,
  \beq
    \label{eq4.9}
        \begin{split}
          \theta^2\phi\eta_1^2 \overline{\hat{q}} \hat{G}^1
           =&\theta^2\phi\eta_1^2 \overline{\hat{q}} \cdot \Big [\hat{q}_{ss}+\sum_{j,k=1}^n(g^{jk}\hat{q}_{x_j})_{x_k}\Big]
           \\
           =&(\theta^2\phi\eta_1^2 \overline{\hat{q}}\hat{q}_s)_s
           -\theta^2\phi\eta_1^2 |\hat{q}_s|^2
           -(\theta^2\phi\eta_1^2)_s\overline{\hat{q}}\hat{q}_s 
           +\sum_{k=1}^n \Big [\theta^2\phi\eta_1^2\sum_{j=1}^ng^{jk}\overline{\hat{q}}\hat{q}_{x_j} \Big]_{x_k}\\
           &-\theta^2\phi\eta_1^2\sum_{j,k=1}^ng^{jk}\hat{q}_{x_j}\overline{\hat{q}}_{x_k}
           -\sum_{j,k=1}^ng^{jk}(\theta^2\phi\eta_1^2)_{x_j}\overline{\hat{q}}\hat{q}_{x_j}.
        \end{split}
    \eeq
It reduces equivalently to 
  \beq
    \label{eq4.9'}
        \begin{split}
            &   \theta^2\phi\eta_1^2 
             \cdot \Big[ |\hat{q}_s|^2  +\sum_{j,k=1}^ng^{jk}\hat{q}_{x_j}\overline{\hat{q}}_{x_k} \Big]
           \\
           =&  \theta^2\phi\eta_1^2 \overline{\hat{q}} \hat{G}^1
            +(\theta^2\phi\eta_1^2 \overline{\hat{q}}\hat{q}_s)_s
           -(\theta^2\phi\eta_1^2)_s\overline{\hat{q}}\hat{q}_s 
           +\sum_{k=1}^n \Big [\theta^2\phi\eta_1^2\sum_{j=1}^ng^{jk}\overline{\hat{q}}\hat{q}_{x_j} \Big]_{x_k}\\
           &
           -\sum_{j,k=1}^ng^{jk}(\theta^2\phi\eta_1^2)_{x_j}\overline{\hat{q}}\hat{q}_{x_j}.
        \end{split}
    \eeq    
   Integrating \eqref{eq4.9'} over~ $(-b,b)\times\Omega$,  we obtain, upon using \eqref{define4.3}, \eqref{principle_coefficient}, \eqref{chapter four positive} and Cauchy-Schwarz inequality, that %
  \be \small
    \label{eq4.12}
        \begin{split}
            \int_{-b}^b\int_{\omega_0} \theta^2\phi(a |\nabla \hat{q}|^2+|\hat{q}_s|^2)dxds 
            &\leq   \int_{-b}^b\int_{\Omega} \theta^2\phi \eta_1^2(a |\nabla \hat{q}|^2+|\hat{q}_s|^2)dxds 
            \\
            &\leq \frac{C}{\lambda\mu^2} \int_{-b}^b\int_{\Omega}\theta^2|\hat G^1|^2dxds
            +C\lambda^2\mu^2\int_{-b}^b\int_{\omega_1}\theta^2 \phi^3 |\hat{q}|^2dxds.
        \end{split}
    \ee
The factor $\frac{C}{\lambda \mu^2}$ is important for the estimate of $I_2$, see \eqref{I2-bound}.

{\bf $\circ $ Step 4.2.}  We estimate a weighted  energy for $\hat q$, i.e.,  $ \int_{-b}^b\int_{\omega_1}\theta^2\phi^3|\hat{q}|^2dxds$. 

 By definition of $\eta_2$, it suffices to find an up bound for 
$ \int_{-b}^b\int_{\Omega}\theta^2\phi^3\eta^2_2|\hat{q}|^2dxds$.
 We claim that there exists a constant $C>0$, such that
    \be
    \label{eq4.13}
        \int_{-b}^b\int_{\Omega}\theta^2\phi^3\eta^2_2|\hat{q}|^2dxds
        \leq C\int_{-b}^b\int_{\Omega}\theta^2\phi^3\eta^2_2|\hat{q}_s|^2dxds.
    \ee
    
Indeed,  in view of  \eqref{eq5.4},  we have  
\beq 
|\hat{q}|=|\varphi z^0|,\ |\hat{q}_s|=|\varphi_s+i\gamma \varphi||z^0|,
\quad  \forall  (s,x)\in [-b,b]\times \Omega. 
\eeq
Notice that 
\beq
 |\varphi_s+i\gamma \varphi|=|\gamma||\varphi|\geq C|\varphi|,~\quad  \forall s \in (-b_0.b_0), \gamma\in \mathcal{O}_C,
 \eeq
 then 
 \be\lab{varphi11}\begin{split}
  \int_{-b_0}^{b_0}\int_{\Omega}\theta^2\phi^3\eta^2_2|\hat{q}|^2dxds
  \leq 
  C \int_{-b_0}^{b_0}\int_{\Omega}\theta^2\phi^3\eta^2_2|\hat{q}_s|^2dxds  
 \leq C \int_{-b}^b\int_{\Omega}\theta^2\phi^3\eta^2_2|\hat{q}_s|^2dxds.
 \end{split}
 \ee 
 Recalling \eqref{the2} and \eqref{the1}, 
we can obtain from the choice of $\varphi$, $b$, $b_0$ in \e{eq4.3} and \e{b-b0},  that
\be \lab{varphi12} 
\begin{split}
\int_{(-b,-b_0)\bigcup(b_0,b)}\int_{\Omega}\theta^2\phi^3\eta^2_2{ |\hat{q}|^2}dxds
= &
\int_{(-b,-b_0)\bigcup(b_0,b)}\int_{\Omega}\theta^2\phi^3\varphi^2\eta^2_2|z^0|^2dxds \\
=&
\int_{\Omega}\eta^2_2|z^0|^2\int_{(-b,-b_0)\bigcup(b_0,b)}\theta^2\phi^3\varphi^2dsdx \\
\leq &
\int_{\Omega}\eta^2_2|z^0|^2\int_{(-b,-b_0)\bigcup(b_0,b)}e^{2\lambda(1+e^\mu)} e^{3\mu \psi(\pm 1,x)}dsdx \\
\leq &
\int_{\Omega}\eta^2_2|z^0|^2\int_{-1}^1e^{2\lambda(2+e^\mu)} e^{3\mu \psi(\pm 1,x)}dsdx \\
\leq &
\int_{\Omega}\eta^2_2|z^0|^2\int_{-1}^1\theta^2\phi^3\varphi^2dsdx\\
\leq &
\int_{-b_0}^{b_0}\int_{\Omega}\theta^2\phi^3\eta^2_2|\hat{q}|^2dxds.
 \end{split}
\ee
Combing \eqref{varphi11} and \eqref{varphi12}  yields  \eqref{eq4.13} immediately. %

{\bf $\circ $ Step 4.3.}     
We estimate ~$\int_{-b}^b\int_{\Omega}\theta^2\phi^3\eta^2_2|\hat{q}_s|^2dxds$ by using the coupling relation in \eqref{define4.3}.  

Multiplying $\hat{p}$-equation in  \eqref{define4.3} by $i\theta^2\phi^3\eta_2^2\overline{\hat{q}}_s$, and arranging the terms,  we get 
  \be\displaystyle 
    \label{eq4.17}
    \begin{split}
        \alpha(x)\theta^2\phi^3\eta_2^2|\hat{q}_s|^2=&i\theta^2\phi^3\eta_2^2 \overline{\hat{q}}_s \Big [\hat{p}_{ss}+\sum_{j,k=1}^n(g^{jk}\hat{p}_{x_j})_{x_k} \Big]\\
        &-i\theta^2\phi^3\eta_2^2 \overline{\hat{q}}_s\, [\hat{G}_1^0+\varphi G^0-i\beta(x)\hat{p}_s],
    \end{split}
    \ee
where $\hat G^0_1$ is given in \e{3.36} and it only concerns the terms with derivates of $\varphi$.    

Direct computations yield that 
{      \be\displaystyle \small
    \label{eq4.15}
        \begin{split}
            &\theta^2\phi^3\eta_2^2\overline{\hat{q}}_s \Big [\hat{p}_{ss}+\sum_{j,k=1}^n(g^{jk}\hat{p}_{x_j})_{x_k} \Big] \\
          =&(\theta^2\phi^3\eta_2^2 \overline{\hat{q}}_{s}\hat{p}_s)_s-(\theta^2\phi^3\eta_2^2)_s\overline{\hat{q}}_{s}\hat{p}_s
         	-\theta^2\phi^3\eta_2^2 \hat{p}_s  \Big[\overline{\hat{q}}_{ss} +\sum_{j,k} (g^{jk}\overline{\hat{q}}_{x_k})_{x_j} \Big]\\
          &+\sum_{j,k}(\theta^2\phi^3\eta_2^2\overline{\hat{q}}_{s}g^{jk}\hat{p}_{x_j})_{x_k}-\sum_{j,k}(\theta^2\phi^3\eta_2^2)_{x_k}\overline{\hat{q}}_{s}g^{jk}\hat{p}_{x_j}
          -\sum_{j,k}(\theta^2\phi^3\eta_2^2\overline{\hat{q}}_{x_k}g^{jk}\hat{p}_{x_j})_s \\
          &+\sum_{j,k}(\theta^2\phi^3\eta_2^2)_s\overline{\hat{q}}_{x_k}g^{jk}\hat{p}_{x_j}+
         \sum_{j,k}(\theta^2\phi^3\eta_2^2 g^{jk}\overline{\hat{q}}_{x_k}\hat{p}_{s})_{x_j}  -\sum_{j,k}(\theta^2\phi^3\eta_2^2)_{x_j} g^{jk}\overline{\hat{q}}_{x_k}\hat{p}_{s} \\
        \end{split}
    \ee}
and
  \be\displaystyle\small
  \label{new1}
        \begin{split}
  \sum_{j,k}(\theta^2\phi^3\eta_2^2)_{x_j} g^{jk}\overline{\hat{q}}_{x_k}\hat{p}_{s}=&\sum_{j,k}[(\theta^2\phi^3\eta_2^2)_{x_j} g^{jk}\overline{\hat{q}}\hat{p}_{s}]_{x_k}-\sum_{j,k}[(\theta^2\phi^3\eta_2^2)_{x_j} g^{jk}\overline{\hat{q}}\hat{p}_{x_k}]_{s}\\
  &-\sum_{j,k}[(\theta^2\phi^3\eta_2^2)_{x_j} g^{jk}]_{x_k}\overline{\hat{q}}\hat{p}_{s}+\sum_{j,k}[(\theta^2\phi^3\eta_2^2)_{x_j} g^{jk}]_{s}\overline{\hat{q}}\hat{p}_{x_k}\\&+\sum_{j,k}(\theta^2\phi^3\eta_2^2)_{x_j} g^{jk}\overline{\hat{q}}_{s}\hat{p}_{x_k}.
     \end{split}
  \ee
 We integrate  \eqref{eq4.17}  on~$(-b,b)\times \Omega$. Then we derive from  \eqref{eq4.15}-\eqref{new1} and particularly the facts that  $\hat p_s = \hat q_s=0 $ on the boundary of $(-b,b) \times \Omega$ and  
 $\nabla \hat p(\pm b,\cdot) = \nabla \hat q_s(\pm b,\cdot) =0$,  that 
    \be
    \label{eq4.18} 
        \begin{split} 
        & \int_{-b}^b\int_{\Omega} \alpha(x) \theta^2\phi^3\eta_2^2 |\hat{q}_s|^2dxds \\
         	\leq & \frac{\delta}{2} \int_{-b}^b\int_{\Omega}\theta^2\phi^3\eta_2^2 |\hat{q}_s|^2dxds  
	        +C  \int_{-b}^b\int_{\omega_2}\theta^2\phi^3\eta_2^2 (|\hat{G}_1^0|^2+| G^0|^2 + |\hat{G}_1^1|^2+| G^1|^2)dxds\\
		&+C \lambda^2\mu^2\int_{-b}^b\int_{\omega_3}\theta^2\phi^5( |\hat{p}|^2+|\hat{p}_s|^2+|\nabla \hat{p}|^2) dxds\\ 
   &+\Big |\int_{-b}^b\int_{\Omega}\sum_{j,k}[(\theta^2\phi^3\eta_2^2)_{x_j} g^{jk}]_{s}\overline{\hat{q}}\hat{p}_{x_k}dxds \Big|
   +\Big|\int_{-b}^b\int_{\Omega}\sum_{j,k}[(\theta^2\phi^3\eta_2^2)_{x_j} g^{jk}]_{x_k}\overline{\hat{q}}\hat{p}_{s}dxds\Big|.
        \end{split}
    \ee
 Since $\alpha(x) \geq  \delta >0$ on $\omega_{\delta}$, the first term on the right hand side of \eqref{eq4.18}  can be absorbed by the left hand side. Moreover,
  by direct computations,  we have the following two facts:
  \be
  \begin{split}
  [(\theta^2\phi^3\eta_2^2)_{x_j} g^{jk}]_{x_k}=&(\theta^2\phi^3\eta_2^2)_{x_j,x_k}g^{jk}+(\theta^2\phi^3\eta_2^2)_{x_j} g^{jk}_{x_k} \\
  =& g^{jk}\big[ 2\theta_{x_j,x_k}\theta\phi^3\eta_2^2+2\theta_{x_j}\theta_{x_k}\phi^3\eta_2^2+6\theta_{x_j}\theta\phi^2\phi_{x_k}\eta_2^2\\
  &+4\theta_{x_j,x_k}\theta\phi^3(\eta_2)_{x_k}\eta_2+6\theta_{x_k}\theta\phi_{x_j}\phi^2\eta_2^2+3\theta^2\phi_{x_j,x_k}\phi^2\eta_2^2\\
  &+6\theta^2\phi_{x_j}\phi\phi_{x_k}\eta_2^2+6\theta^2\phi_{x_j}\phi^2\eta_2(\eta_2)_{x_k}+2(\eta_2)_{x_j,x_k}\eta_2\theta^2\phi^3\\
  &+2(\eta_2)_{x_j}(\eta_2)_{x_k}\theta^2\phi^3+4(\eta_2)_{x_j}\eta_2\theta_{x_k}\theta\phi^3+6(\eta_2)_{x_j}\eta_2\theta^2\phi^2\phi_{x_k}\big]\\
  &+[2\theta_{x_j}\theta\phi^3\eta_2^2+3\theta^2\phi_{x_j}\phi^2\eta_2^2+2(\eta_2)_{x_j}\eta_2\theta^2\phi^3 ]g^{jk}_{x_k}
  \end{split}
  \ee  
  and 
  \be
  \begin{split}
   [(\theta^2\phi^3\eta_2^2)_{x_j} g^{jk}]_{s}=&g^{jk}[2\theta_{x_j,s}\theta\phi^3\eta_2^2+2\theta_{x_j}\theta_s\phi^3\eta_2^2+6\theta_{x_j}\theta\phi_s\phi^2\eta_2^2 \\
   &+6\theta_s\theta\phi_{x_j}\phi^2\eta_2^2+3\theta^2\phi_{x_j,s}\phi^2\eta_2^2+6\theta^2\phi_{x_j}\phi_s\phi\eta_2^2  \\
   &+2(\eta_2)_{x_j}\eta_2\theta_s\theta\phi^3+6(\eta_2)_{x_j}\eta_2\theta^2\phi^2\phi_s].
   \end{split}
  \ee
  By definition of $\eta_2, \theta, \phi$, in particular \eqref{eta2}, \eqref{eq4.7} and $g^{jk}(.)\in C^1(\overline{\Omega};\mathbb{R})$, one can obtain the following estimate:
  \be\label{3.52}
\big| [(\theta^2\phi^3\eta_2^2)_{x_j} g^{jk}]_{s}\big|
+   \big| [(\theta^2\phi^3\eta_2^2)_{x_j} g^{jk}]_{x_k}\big|\leq C\lambda^2\mu^2\eta_2\theta^2  \phi^5.
   \ee
Thus by above estimate \eqref{3.52} and Cauchy-Schwarz inequality, we get for all $\varepsilon>0$ that 
  \beq
    \label{3.62} 
        \begin{split} 
       & \Big |\int_{-b}^b\int_{\Omega}\sum_{j,k}[(\theta^2\phi^3\eta_2^2)_{x_j} g^{jk}]_{s}\overline{\hat{q}}\hat{p}_{x_k}dxds \Big|
        +\Big|\int_{-b}^b\int_{\Omega}\sum_{j,k}[(\theta^2\phi^3\eta_2^2)_{x_j} g^{jk}]_{x_k}\overline{\hat{q}}\hat{p}_{s}dxds\Big|
         \\
        	\leq &  \varepsilon  \int_{-b}^b\int_{\Omega}\theta^2\phi^3\eta_2^2 |\hat{q}|^2dxds  
	+ \frac{C \lambda^4\mu^4}{\varepsilon} \int_{-b}^b\int_{\omega_3}\theta^2\phi^7(|\hat{p}_s|^2+|\nabla \hat p|^2) dxds,
        \end{split}
    \eeq
 Therefore we obtain  from \eqref{eq4.13}  and by choosing $\varepsilon>0$ small that 
 \be
    \label{new2}
        \begin{split}
      	   \int_{-b}^b\int_{\Omega}\theta^2\phi^3\eta^2_2|\hat{q}_s|^2dxds
	       \leq & C  \int_{-b}^b\int_{\Omega}\theta^2\phi^3\eta_2^2(|\hat{G}_1^0|^2+| G^0|^2 + |\hat{G}_1^1|^2+| G^1|^2 )dxds\\
		&+C\lambda^4\mu^4\int_{-b}^b\int_{\omega_3}\theta^2\phi^{7}( |\hat{p}|^2+|\hat{p}_s|^2+|\nabla \hat{p}|^2) dxds,
        \end{split}
    \ee
 \be
    \label{new3}
        \begin{split} 
         \int_{-b}^b\int_{\Omega}  \theta^2\phi^3\eta_2^2 |\hat{q}|^2dxds
	   \leq & C  \int_{-b}^b\int_{\Omega}\theta^2\phi^3\eta_2^2(|\hat{G}_1^0|^2+| G^0|^2 + |\hat{G}_1^1|^2+| G^1|^2 )dxds\\
		&+C\lambda^4\mu^4\int_{-b}^b\int_{\omega_3}\theta^2\phi^7( |\hat{p}|^2+|\hat{p}_s|^2+|\nabla \hat{p}|^2) dxds.
        \end{split}
    \ee

{\bf $\circ $  Step 4.4.}  We estimate  $ \int_{-b}^b\int_{\omega_3} \theta^2\phi^7( |\nabla \hat{p}|^2+|\hat{p}_s|^2)dxds$.

Similar  to \eqref{eq4.12}, we can derive 
    \be\displaystyle\small
    \label{eq4.10} 
        \begin{split} 
            \int_{-b}^b\int_{\omega_3} \theta^2\phi^7(|\nabla \hat{p}|^2+|\hat{p}_s|^2)dxds
            \leq \frac{C}{\lambda^7\mu^8} \int_{-b}^b\int_{\Omega}\theta^2 |\hat G^0|^2dxds 
            +C \lambda^7\mu^{8}\int_{-b}^b\int_{\omega_\delta}\theta^2\phi^{14}|\hat{p}|^2dxds .
        \end{split}
    \ee
The factor $\frac{C}{\lambda^7 \mu^8}$ is important for the estimate of $I_2$, see \eqref{I2-bound}.    

{\bf $\circ $  Step 4.5.}  We summarize the estimate of $I_2$ from above.

Applying \eqref{eq4.12}, \eqref{new2},\eqref{new3}, \eqref{eq4.10}  for $I_{23},I_{24},I_{21}$ in \eqref{I2}, we end up with
\be 
    \label{new5}
        \begin{split} 
      I_{2} = & I_{21}+I_{22}+I_{23}+I_{24}\\
    	     \leq & C \lambda^3 \mu^4 \int_{-b}^b\int_{\Omega}
	   	\theta^2\phi^3\eta_2^2 (|\hat{G}_1^0|^2+ |\hat{G}_1^1|^2 )dxds
	     +C  \int_{-b}^b\int_{\Omega}\theta^2(|\hat{G}^0|^2 +|\hat G^1|^2 )dxds\\
	       	     & +C \lambda^3 \mu^4 \int_{-b}^b\int_{\Omega}
	   	\theta^2\phi^3\eta_2^2 (| G^0|^2 +| G^1|^2 )dxds
		+C\lambda^{14} \mu^{16} \int_{-b}^b\int_{\omega_{\delta}}\theta^2\phi^{14}|\hat{p}|^2 dxds.
        \end{split}
    \ee
Similar to the estimate of $I_{11}$ in \eqref{I11}, we have  
 \beq   \displaystyle
      \begin{split}
             C \lambda^3 \mu^4 \int_{-b}^b\int_{\Omega}
	   	\theta^2\phi^3\eta_2^2 (|\hat{G}_1^0|^2+ |\hat{G}_1^1|^2 )dxds
              \leq \lambda^3 e^{2\lambda(1+e^{\mu})}C(\mu) (\|p\|^2_{H^1(X)} +\|q\|^2_{H^1(X)}) .
  \end{split}
    \eeq
Note that $C  \int_{-b}^b\int_{\Omega}\theta^2(|\hat{G}^0|^2 +|\hat G^1|^2 )dxds$ is the same as $I_1$ and the estimate has been given  by \eqref{I1-bound}.  Therefore, we have 
\be 
    \label{I2-bound}
        \begin{split} 
      I_{2} \leq & CI_1+ \lambda^3 e^{2\lambda(1+e^{\mu})}C(\mu) (\|p\|^2_{H^1(X)} +\|q\|^2_{H^1(X)}) \\
	       	     & +C \lambda^3 \mu^4 \int_{-b}^b\int_{\Omega}
	   	\theta^2\phi^3\eta_2^2 (| G^0|^2 +| G^1|^2 )dxds
	+C\lambda^{14} \mu^{16} \int_{-b}^b\int_{\omega_{\delta}}\theta^2\phi^{14}|\hat{p}|^2 dxds.
        \end{split}
    \ee

{\bf $\circ $ Step 5.}  We derive the estimate \eqref{eq4.1} based on Step 1-4. 

Letting $\lambda$ be large enough,  we finally obtain from \eqref{3.35}, \eqref{I0-bound}, \eqref{I1-bound}  and  \eqref{I2-bound},  the desired interpolation estimate \eqref{eq4.1}. %

%
%
%

{\bf $\circ $ Step 6.} We provide an example of the cut-off function $\eta_2 \in C_0^3(\omega_2)$ such that 
\eqref{eq4.8}-\eqref{eta2} hold.
 
 Without loss of generality, we may assume  $\omega_1\subset B(0,\frac{r}{2})\subset B(0,r)\subset \omega_2$. 
Set
  \beq
    \eta_2(x) \triangleq
  \begin{cases}
1,~& \text{if\ } |x|< \frac{r}{2}. \\
\frac{128}{81r^{10}}[r^2-|x|^2]^3[16|x|^4- 2r^2|x|^2+r^4], ~& \text{if\ }  \frac{r}{2} \leq |x| \leq r,\\
0,~& \text{if\ }  |x|> r. \\  
\end{cases}
\eeq
It is easy to check that $\eta_2 \in C_0^3(\omega_2)$ and 
\beq
\lim_{|x|\rightarrow r} \frac{|(\eta_2)_{x_k}(\eta_2)_{x_j}|}{\eta_2}=0, 
\quad 
\lim_{|x|\rightarrow \frac{r}{2}} \frac{|(\eta_2)_{x_k}(\eta_2)_{x_j}|}{\eta_2}<\infty.
\eeq
Then it follows the property \eqref{eta2}.

This finally concludes the proof of Lemma \ref{lem-key}. \hfill$\Box$

\section{An example with $\omega_{\alpha}\cap \omega_{\beta} =\emptyset$}\lab{sect:ex}

In this section, we show, by an example, that the geometric assumption $\omega_{\alpha}\cap \omega_{\beta} \neq \emptyset$ is necessary in general for 
asymptotic stability of the coupled system \eqref{wave}.

More precisely, we consider the following 1-D coupled wave equations
\be\lab{ex:cou}
\begin{cases}
y_{tt} -y_{xx}+\alpha(x)z_t +\beta(x) y_t=0,  \quad &\text{in} ~(0,+\infty)\times (0,2\pi),\\
z_{tt}- z_{xx}-\alpha(x)y_t=0, \quad &\text{in} ~(0,+\infty)\times (0,2\pi),\\
u(-\pi)=u(\pi)=v(-\pi)=v(\pi)=0, \quad &\text{in} ~(0,+\infty),\\
 (y(0),y_t(0),z(0),z_t(0))=(y^0,y^1,z^0,z^1), \quad &\text{in}~ (0,2\pi),
\end{cases}
\ee
where the coefficients $\alpha$ and $\beta$ are given by 
\beq
\lab{al}
\alpha(x)=\left\{\begin{split}
0, \quad&\text{if}~ x\in (\pi,2\pi), \\
\f{24}{5}, \quad&\text{if}~ x\in [0,\pi),
\end{split}\right.  \quad 
\beta(x)=\begin{cases}
1, \quad \text{if}~ x\in (\pi,2\pi),\\
0, \quad \text{if}~ x\in [0,\pi)
\end{cases}
\eeq
such that  $\omega_{\alpha} \cap \omega_{\beta}  = \emptyset.$

%
%

Let the initial data be the following 
 \be\lab{ex:in1}
   y^0(x)=0, \quad \text{if}~ x\in (0,2\pi), \quad     y^1(x)=
       \begin{cases}
        7sin x- sin 7x,&\text{if}~ x\in (0,\pi),\\
        0,&\text{if}~ x\in [\pi,2\pi),\\
       \end{cases},\quad   
      \ee
and 
 \be\lab{ex:in2}
     z^0(x)=   \begin{cases}
        -7sin x- sin 7x,&\text{if}~ x\in (0,\pi),\\
        -\frac{14}{5}sin 5x,&\text{if}~ x\in [\pi,2\pi),\\
       \end{cases},\quad    
      z^1(x)=0 ,\quad \text{if}~ x\in (0,2\pi),
      \ee
 Clearly, one can see that $(y^0,y^1,z^0,z^1)\in \mathcal{D}(\mathcal{A})$, the solution of the system \eqref{ex:cou}-\eqref{ex:in1} can be obtained 

    \be\lab{cont}
  y(t,x)=
   \begin{cases}
      sin (5t) (7sin x- sin 7x)  ,&\text{if}~ x\in (0,\pi),\\
     0,&\text{if}~ x\in [\pi,2\pi),\\
    \end{cases}
    \ee
   and 
   \be
   z(t,x)=   \begin{cases}
   -cos (5t) (7sin x+sin (7x))  ,&\text{if}~ x\in (0,\pi),\\
     -\frac{14}{5}cos (5t) sin (5x),&\text{if}~ x\in [\pi,2\pi),\\
    \end{cases}
   \ee
In contrast to \eqref{eng-decay} in Theorem \ref{mainthm}, the energy of the system \eqref{ex:cou} is conserved :
  \beq  \f{d}{dt} \mathbb{E}(y,z)(t)=-\int_{-\pi}^{\pi} \beta(x)|y_t|^2 dx =0, \quad  \forall t \geq 0,    \eeq
therefore the system is not asymptotically stable.  Furthermore,  decay rate estimate \eqref{logstable} is invalid for some initial data $(y^0,y^1,z^0,z^1)\in \mathcal{D(A)}$. 
%

\medskip

\textbf{Acknowledgements:}
The first author is partially supported by the fund of the Shanghai Key Laboratory for Contemporary Applied Mathematics at Fudan University (No. 74120-42080001).
The second author is partially supported by the Natural Science Foundation of China (No. 11971119). 

\def\cprime{$'$}

\end{CJK*}

\end{document}